\documentclass[9pt,shortpaper,web,twoside]{IEEEtran}
\newtheorem{definition}{\bf Definition}
 
  \newtheorem{Corollary}{\bf Corollary}
  \newtheorem{thm}{\bf Theorem}

 \newtheorem{lemma}{\bf Lemma}
  \newtheorem{prop}{\bf Proposition}
\usepackage{amsfonts,amssymb,amsmath,graphicx}  
\usepackage{mathtools}
\usepackage{bbm}
\usepackage{breqn}
\usepackage{subcaption}
\usepackage{diagbox}
\usepackage{cite}
\usepackage{cases}
\usepackage{tabu}
\usepackage{url}
\usepackage{multicol}
\usepackage[bottom]{footmisc}
\usepackage{environ}
\usepackage{tikz} 
 \usepackage{pgfplots}
\pgfplotsset{compat=newest}
\usetikzlibrary{plotmarks}
\usetikzlibrary{arrows.meta}
\usepgfplotslibrary{patchplots}
\usetikzlibrary{arrows,decorations.markings}
\usetikzlibrary{decorations.pathreplacing}
\usetikzlibrary{automata,arrows,positioning,calc}
\usepackage{grffile}
\usepackage{enumerate}
\date{}
\title{On Minimizing Total Discounted Cost in MDPs\\
Subject to Reachability Constraints}
\author{Yagiz Savas, Christos K. Verginis, Michael Hibbard, and Ufuk Topcu \thanks{All authors are with the University of Texas at Austin, Austin, TX 78705 USA. E-mail:\{yagiz.savas, mhibbard, utopcu\}@utexas.edu,  christos.verginis@austin.utexas.edu.} \thanks{This work was supported in part by the grants ARL W911NF-17-2-0181, AFRL FA9550-19-1-0169, and DARPA D19AP00004.}}
\begin{document}
\maketitle

\begin{abstract}
   We study the synthesis of a policy in a Markov decision process (MDP) following which an agent reaches a target state in the MDP while minimizing its total discounted cost. The problem combines a reachability criterion with a discounted cost criterion and naturally expresses the completion of a task with probabilistic guarantees and optimal transient performance. We first establish that an optimal policy for the considered formulation may not exist but that there always exists a near-optimal stationary policy. We additionally provide a necessary and sufficient condition for the existence of an optimal policy. We then restrict our attention to stationary deterministic policies and show that the decision problem associated with the synthesis of an optimal stationary deterministic policy is NP-complete. Finally, we provide an exact algorithm based on mixed-integer linear programming and propose an efficient approximation algorithm based on linear programming for the synthesis of an optimal stationary deterministic policy. 
\end{abstract}
\section{Introduction}\label{introduction_sec}
Markov decision processes (MDPs) provide a framework to model the behavior of an agent, e.g., humans or autonomous robots, operating under uncertainty  \cite{puterman2014markov}. A typical objective in an MDP is to synthesize a policy under which
the agent 
reaches a set of target states with \textit{maximum probability} \cite{de1999computing}. Such a reachability objective may express, e.g., the completion of a surveillance mission in robotics applications \cite{ding2014optimal}, the achievement of a drug concentration level in blood in healthcare applications \cite{hu1996comparison,schaefer2005modeling}, and the alignment of a portfolio with investors' preferences in finance applications \cite{pola2011stochastic}. 

In many planning problems, achieving a desired transient behavior is as important as the task completion.
A widely used performance criterion in MDPs is the total \textit{discounted} cost 
accumulated by the agent along its trajectories \cite{puterman2014markov,altman1999constrained}. By associating the agent's actions with \textit{nonnegative} costs and discounting the costs incurred in the future, such a criterion naturally expresses the agent's short- and long-term 
considerations. 
Discounting may represent, e.g., the importance of early detection of an intruder in robotic applications \cite{kiennert2018survey} or the agent's opportunity cost in healthcare \cite{schaefer2005modeling} and finance \cite{torgerson1999discounting} applications. 


In this paper, we present a comprehensive analysis for the problem of synthesizing a policy under which an agent reaches a desired set of target states with maximum probability while minimizing its total discounted cost. Such a policy ensures the completion of a task with probabilistic guarantees. Moreover, when the task can be completed at a minimum total cost by following multiple trajectories, the synthesized policy allows the agent to decide on the ordering of the events happening until completion thanks to discounting.

In the literature, several problem formulations have been proposed to synthesize policies satisfying multiple criteria. Extensively studied problem formulations include the so-called constrained MDP problems \cite{altman1999constrained}, stochastic shortest path (SSP) problems \cite{bertsekas1991analysis}, and multi-objective model checking problems \cite{etessami2007multi}. The constrained MDP problem associates
the agent's actions with multiple costs
that are discounted with either the same \cite{altman1999constrained,chatterjee2006markov} or different \cite{dolgov2005stationary,feinberg1999constrained,feinberg2000constrained,chen2004dynamic} discount factors. In general, probabilistic reachability objectives cannot be expressed as total discounted criteria unless one makes restrictive assumptions on the MDP structure \cite{feinberg2019reduction}. 
In the SSP \cite{bertsekas1991analysis,Bertsekas,teichteil2012stochastic} and multi-objective model checking problems  \cite{etessami2007multi,hartmanns2018multi,delgrange2020simple}, one aims to synthesize a policy that satisfies certain reachability constraints while minimizing the agent's total \textit{undiscounted} cost. Due to the lack of discounting, however, the agent cannot adjust the importance of its short- and long-term considerations. Overall, the existing methods for multi-objective planning fail when one needs to synthesize policies that minimize the total discounted cost while satisfying a task with probabilistic guarantees. 


It is known that an optimal policy always exists for the previously described problem formulations \cite{altman1999constrained,feinberg1992constrained,bertsekas1991analysis,etessami2007multi}. However, to the best of our knowledge, the existence of optimal solutions for the formulation studied in this paper is an open problem.
Our first contribution is establishing that, in an MDP, a policy that minimizes 
the total discounted cost among the ones that maximize 
the probability of reaching a target state may not exist. This result illustrates a fundamental difference of the problem considered in this paper from the formulations considered in the literature. 

Since optimal policies do not exist for the general case, it is critical to verify whether there exists an optimal policy for a given problem instance. As the second contribution, we present an efficiently verifiable 
necessary and sufficient condition for the existence of optimal policies in a given problem instance. When there are no optimal policies, one typically searches for near-optimal policies. As the third contribution, we show that, for any positive constant $\epsilon$, 
there exists an $\epsilon$-optimal \textit{stationary} policy that can be synthesized efficiently.

In many applications, it is desirable to generate
stationary deterministic policies
due to their low computational requirements. 
It is known \cite{feinberg1999constrained} that, in general, the synthesis of such policies is NP-hard for constrained MDPs in which multiple costs are discounted with the same discount factor. Since we consider a problem involving a probabilistic constraint, however, the existing complexity results on constrained MDPs do not apply to 
the problem considered in this paper. As the fourth contribution, we establish that it is NP-complete to decide whether there exists a \textit{stationary deterministic} policy that maximizes the probability of reaching a set of target states while attaining a total discounted cost below a desired threshold. 
Motivated by this complexity result, we synthesize an optimal stationary deterministic policy by formulating a mixed-integer linear program (MILP) which is an extension of the algorithms developed in \cite{dolgov2005stationary,kalagarla2020synthesis}. For small problem instances, one can compute optimal solutions to MILPs using off-the-shelf solvers, e.g., \cite{gurobi}. However, for large problem instances, MILP formulations become 
intractable. To remedy this limitation, our fifth contribution is the development of an approximation algorithm based on linear programming that efficiently synthesizes stationary deterministic policies with theoretical suboptimality guarantees. 

In numerical simulations, we present an application of the studied formulation to motion planning. Specifically, we consider an agent that aims to deliver a package to a certain location while minimizing the risk of being attacked by adversaries in the environment. We illustrate that the proposed approximation algorithm generates agent trajectories that guarantee reachability to the desired location while visiting the minimum number of risky regions in the environment.


\section{Preliminaries}\label{prelim_sec}
 \noindent\textbf{Notation:} The sets of natural and real numbers are denoted by $\mathbb{N}$ and $\mathbb{R}$, respectively, while the set of nonnegative reals is denoted by $\mathbb{R}_{\geq 0}$. Finally, $|S|$ denotes the cardinality of a set $S$.
{\setlength{\parindent}{0cm}
\begin{definition}
A \textit{Markov decision process} (MDP) is a tuple $\mathcal{M}$$:=$$(S, s_1, \mathcal{A}, \mathcal{P})$ where $S$ is a finite set of states, $s_1$$\in$$S$ is an initial state, $\mathcal{A}$ is a finite set of actions, and $ \mathcal{P}$$:$$S$$\times$$ \mathcal{A}$$\times$$S$$\rightarrow$$[0,1]$ is a transition function such that $\sum_{s'\in S}\mathcal{P}(s,a,s')$$=$$1$ for all $s$$\in$$S$ and $a$$\in$$\mathcal{A}(s)$, where $\mathcal{A}(s)$$\subseteq$$\mathcal{A}$ denotes the set of available actions in $s$$\in$$S$.
\end{definition}}

We denote the transition probability $ \mathcal{P}(s,a,s')$ by $ \mathcal{P}_{s,a,s'}$. A state $s$$\in$$S$ is \textit{absorbing} if $\mathcal{P}_{s,a,s}$$=$$1$ for all $a$$\in$$\mathcal{A}(s)$.
{\setlength{\parindent}{0cm}
\noindent \begin{definition}
For an MDP $\mathcal{M}$, a \textit{policy} $\pi$$:=$$(d_1,d_2, d_3, \ldots)$ is a sequence where, for each $t$$\in$$\mathbb{N}$, $d_t$$:$$S$$\times$$\mathcal{A}$$\rightarrow$$[0,1]$ is a mapping such that $\sum_{a\in \mathcal{A}(s)}d_t(s,a)$$=$$1$ for all $s$$\in$$S$. A \textit{stationary} policy is a policy of the form $\pi$$=$$(d_1,d_1, d_1, \ldots)$. A stationary \textit{deterministic} policy is a stationary policy where, for each $s$$\in$$S$, $d_1(s,a)$$=$$1$ for some $a$$\in$$\mathcal{A}(s)$. We denote the set of all policies, all stationary policies, and all stationary deterministic policies by $\Pi(\mathcal{M})$, $\Pi^S(\mathcal{M})$, and $\Pi^{SD}(\mathcal{M})$, respectively.
\end{definition}}

A policy $\pi$$\in$$\Pi(\mathcal{M})$ is traditionally referred to as a Markovian policy \cite{puterman2014markov}. Although it is possible to consider more general policy classes, the consideration of the set $\Pi(\mathcal{M})$ is without loss of generality for the purposes of this paper due to Theorem 5.5.1 in \cite{puterman2014markov}. For notational simplicity, we denote the probability of taking an action $a$$\in$$\mathcal{A}$ in a state $s$$\in$$S$ under a \textit{stationary} policy $\pi$ by $\pi(s,a)$.

A \textit{path} is a sequence $\varrho^{\pi}$$=$$s_1s_2s_3\ldots$ of states generated in $\mathcal{M}$ under $\pi$ which satisfies $\mathcal{P}_{s_t,a_t,s_{t+1}}$$>$$0$ for all $t$$\in$$\mathbb{N}$.
We define the set of all paths in $\mathcal{M}$ under $\pi$ by $Paths^{\pi}_{\mathcal{M}}$ and use the standard probability measure over the set $Paths^{\pi}_{\mathcal{M}}$ \cite{Model_checking}. Let $\varrho^{\pi}[t]$$:=$$s_t$ denote the state visited at the $t$-th step along the path $\varrho^{\pi}$. We define
\begin{align*}
     \text{Pr}_{\mathcal{M}}^{\pi}(Reach[B]):=\text{Pr}\{\varrho^{\pi}\in Paths^{\pi}_{\mathcal{M}}: \exists t\in \mathbb{N}, \varrho^{\pi}[t]\in B \}
\end{align*}
as the probability with which the paths generated in $\mathcal{M}$ under $\pi$ reaches the set $B$$\subseteq$$S$.

\section{Problem Statement}
We consider an agent that aims to reach a set of target states with maximum probability while minimizing the expected total discounted cost it accumulates along its path. Formally, let $c$$:$$S$$\times$$\mathcal{A}$$\rightarrow$$\mathbb{R}_{\geq 0}$ be a \textit{non-negative} cost function, $\beta$$\in$$(0,1)$ be a discount factor, $B$$\subseteq$$S$ be a set of \textit{absorbing} target states, and $J$$:$$\Pi(\mathcal{M})$$\rightarrow$$\mathbb{R}_{\geq 0}$ be a cost-to-go function such that 
\begin{align}
    J(\pi) := \mathbb{E}^{\pi}\Bigg[\sum_{t=1}^{\infty}\beta^{t-1} c(s_t,a_t) \Bigg],
\end{align}
where the expectation is taken over the paths generated under $\pi$. Moreover, let $\Xi(\mathcal{M},B)$$\subseteq$$\Pi(\mathcal{M})$ denote the set of policies under which the paths in $\mathcal{M}$ reach the target set $B$ with maximum probability, i.e.,  $\overline{\pi}$$\in$$\Xi(\mathcal{M},B)$ if and only if
\begin{align}
    \overline{\pi} \in \arg\max_{\pi\in\Pi(\mathcal{M})} \text{Pr}_{\mathcal{M}}^{\pi}(Reach[B]).\label{max_reac_def}
\end{align}
Existence of the maximum in \eqref{max_reac_def} is due to Lemma 10.102 in \cite{Model_checking}.

In this paper, we study the synthesis of a policy $\pi^{\star}$ such that 
\begin{align}\label{main_objective}
  \pi^{\star} \in \arg \inf_{\pi\in\Xi(\mathcal{M},B)}\  J(\pi).
\end{align}


In the following sections, we analyze the problem in 
\eqref{main_objective} and present efficient algorithms to synthesize policies that either exactly or approximately satisfy the condition in \eqref{main_objective}.



\section{Non-Existence of Optimal Policies} \label{non_existence_sec}
We now present a numerical example to illustrate that an optimal policy satisfying the condition in \eqref{main_objective} may not exist.
This example demonstrates the significant difference of the problem considered in this work 
from the traditional constrained MDP problems, e.g., \cite{altman1999constrained,chen2004dynamic, dolgov2005stationary}, for which an optimal policy always exists. 

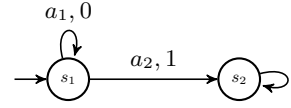
\begin{figure}[t!]
\centering
\begin{tikzpicture}[->, >=stealth', auto, semithick, node distance=2cm]

    \tikzstyle{every state}=[fill=white,draw=black,thick,text=black,scale=0.7]

    \node[state,initial,initial text=] (s_1) {$s_1$};
    \node[state] (s_2) [right=17mm of s_1]  {$s_2$};

\path
(s_1)  edge  [loop above=10]    node{$a_1, 0$}     (s_1)
(s_1)	 edge     node{$a_2, 1$}     (s_2)
(s_2)  edge  [loop right=10]    node{}     (s_2);
\end{tikzpicture}
\caption{An MDP example to illustrate the non-existence of optimal policies. The initial state is $s_1$, and $B$$=$$\{s_2\}$. The tuples $(a,c)$ indicate the action $a$ and the cost $c$.}
\label{fig:optimal_policies}
\end{figure}

Consider the example shown in Fig. \ref{fig:optimal_policies}. The agent starts from the state $s_1$ and aims to reach the state $s_2$ with probability one, i.e., $B$$=$$\{s_2\}$, while minimizing its total discounted cost. Suppose that the agent follows the stationary policy $\pi$$\in$$\Pi^S(\mathcal{M})$ such that $\pi(s_1,a_1)$$=$$1-\delta$ and $\pi(s_1,a_2)$$=$$\delta$
where $\delta$$\in$$[0,1]$. For $\delta$$\in$$(0,1]$, the agent reaches the state $s_2$ with probability one under $\pi$. However, if $\delta$$=$$0$, we have $\pi$$\not \in$$\Xi(\mathcal{M},B)$. Hence, the set $\Xi(\mathcal{M},B)$$\cap$$\Pi^S(\mathcal{M})$ of feasible stationary policies is
\begin{align*}
    \{\pi\in \Pi^S(\mathcal{M}) \ |\  \pi(s_1,a_1)=1-\delta,\ \pi(s_1,a_2)=\delta, \ \delta \in (0,1] \}.
\end{align*}


For a given stationary policy $\pi$$\in$$\Xi(\mathcal{M},B)$$\cap$$\Pi^S(\mathcal{M})$, we have
\begin{align*}
   J(\pi)= \frac{\delta}{1-\beta (1-\delta)}. 
\end{align*}
Note that $\lim_{\delta \rightarrow 0}\frac{\delta}{1-\beta (1-\delta)}$$=$$0$, which implies that
\begin{align*}
   \inf_{\pi\in \Xi(\mathcal{M},B)} J(\pi)=0.
\end{align*}
The \textit{only} policy that attains the infimum is the \textit{stationary} policy $\overline{\pi}$ such that $\overline{\pi}(s_1,a_1)$$=$$1$ and $\overline{\pi}(s_1,a_2)$$=$$0$ since any other policy incurs a non-zero cost by taking the action $a_2$ with non-zero probability. Note that the policy $\overline{\pi}$ is not in the feasible policy space as it reaches the set $B$ with probability zero. Since the infimum is not attainable by any feasible policy, we conclude that there exists no optimal policy.

One may be tempted to think that the existence of a zero cost action $a_1$ is the reason for not having an optimal policy in this example. It can be shown that, even if we assign a positive cost for the action $a_1$, e.g., $c(s_1,a_1)$$=$$0.1$, an optimal policy still does not exist so long as we choose a small discount factor, e.g., $\beta$$<$$0.1$. An optimal policy does not exist in this example because the agent exploits the discounting in the costs and stays in the initial state for as long as possible before reaching the target state. By doing so, the agent ensures that its cost-to-go approaches zero while still satisfying the reachability constraint asymptotically. Such a behavior is specific to the problem in \eqref{main_objective}, which involves a discounted criterion in the objective and a probabilistic criterion in the constraint, and does not arise in existing MDP formulations.

\section{Existence of Near-Optimal Policies}\label{near_opt_sec}
In the previous section, we showed that an optimal policy solving the problem in
\eqref{main_objective} may not exist. Here, we show that there always exists a near-optimal policy which can be synthesized efficiently.
{\setlength{\parindent}{0cm}
\noindent \begin{definition}
For a given constant $\epsilon$$>$$0$, a policy $\overline{\pi}$$\in$$\Xi(\mathcal{M},B)$ is said to be an $\epsilon$-optimal policy for the problem in \eqref{main_objective} if
\begin{align*}
    J(\overline{\pi}) \leq \inf_{\pi\in\Xi(\mathcal{M},B)}J(\pi)+\epsilon.
\end{align*}
\end{definition}}


We establish the existence of $\epsilon$-optimal policies for the problem in \eqref{main_objective} in three steps. First, we \textit{clean up} the MDP $\mathcal{M}$ by removing from the states $s$$\in$$S$ the actions $a$$\in$$\mathcal{A}(s)$ that are guaranteed to yield infeasible policies $\pi$$\in$$\Pi(\mathcal{M}) \backslash \Xi(\mathcal{M},B)$. 
Second, we synthesize a policy on the resulting MDP under which the total discounted cost is minimized in the absence of the reachability constraint. Finally, we perturb the synthesized policy to obtain a policy that  satisfies the reachability constraint and  is $\epsilon$-optimal for the problem in
\eqref{main_objective}.

\subsection{Cleaning up the MDP}\label{cleaned_up_sec}
For a given MDP $\mathcal{M}$, we first partition the set $S$ of states into three disjoint sets. Let $B$$\subseteq$$S$ be the set of target states, and $S_0$$\subseteq$$S$ be the set of states that have zero probability of reaching the states in $B$ under any policy. Finally, we let $S_r$$\coloneqq $$S\backslash (B\cup S_0)$. These sets can be computed efficiently using graph search algorithms \cite{Model_checking}.

Let ${\bf{x}}$$\coloneqq$$(x_s)_{s\in S}$$\in$$\mathbb{R}^{\lvert S\rvert}$ be a vector such that
\begin{align}
    x_s\coloneqq\begin{cases} 1 & \text{if} \ s\in B\\
    0 & \text{if} \ s \in S_0 \\
    \max\limits_{{a\in \mathcal{A}(s)}}\Big\{ \sum_{s'\in S}\mathcal{P}_{s,a,s'} x_{s'} \Big\}& \text{otherwise}.
    \end{cases} \label{x_vector}
\end{align}
It is known  \cite[Chapter 10]{Model_checking} that each element $x_s$ of ${\bf{x}}$ corresponds to the maximum probability of reaching $B$ from $s$$\in$$S$. Moreover, ${\bf{x}}$ can be efficiently computed via linear programming \cite{Model_checking}.

Now, let $S^{\pi}_{\rightarrow}$$\subseteq$$S$ be the set of states that are reachable from the initial state under the policy $\pi$$\in$$\Pi(\mathcal{M})$, i.e.,
\begin{align*}
 S^{\pi}_{\rightarrow} \coloneqq \{s\in S \ |\ \text{Pr}^{\pi}_{\mathcal{M}}(Reach[s]) > 0\}.
\end{align*}
The following result, which is due to Theorem 10.100 in \cite{Model_checking}, characterizes a necessary condition for a policy $\pi$$\in$$\Pi(\mathcal{M})$ to maximize the probability of reaching the target set $B$.
{\setlength{\parindent}{0cm}
\noindent \begin{prop}  \label{prop_1_statement} \cite{Model_checking}
If $\pi$$\in$$\Xi(\mathcal{M},B)$, then, for all $s$$\in$$ S_{\rightarrow}^{\pi}$$\cap$$S_r$,
\begin{align}\label{necc_reach}
    x_s = \sum_{s'\in S}\mathcal{P}_{s,a,s'}x_{s'}
\end{align}
for all $a$$\in$$\mathcal{A}(s)$ satisfying $d_t(s,a)$$>$$0$ for some $t$$\in$$\mathbb{N}$.
\end{prop}}

Equation \eqref{necc_reach} constitutes a necessary (but not sufficient) condition for the feasibility of a policy $\pi$$\in$$\Pi(\mathcal{M})$ for the problem in \eqref{main_objective}. Therefore, without loss of generality, we can remove from the MDP all actions that violate the equality in \eqref{necc_reach}. Specifically, given an MDP $\mathcal{M}$, we obtain the \textit{cleaned-up MDP} $\mathcal{M}'$ by removing all actions $a$$\in$$\mathcal{A}(s)\backslash \mathcal{A}_{\max}(s)$ from each state $s$$\in$$S_r$, where
\begin{align}\label{A_max_def}
    \mathcal{A}_{\max}(s)\coloneqq\Bigg\{a\in\mathcal{A}(s)\  \Big| \ x_s = \sum_{s'\in S}\mathcal{P}_{s,a,s'}x_{s'}\Bigg\}.
\end{align}

\subsection{Minimizing the total cost on the cleaned-up MDP}
On the cleaned-up MDP $\mathcal{M}'$, we synthesize a stationary deterministic policy $\widetilde{\pi}$$\in$$\Pi^{SD}(\mathcal{M}')$ such that
\begin{align}\label{min_unconstrained}
\widetilde{\pi} \in \arg \min_{\pi\in\Pi(\mathcal{M}')}\  J(\pi).
\end{align}
The existence of a stationary deterministic policy $\widetilde{\pi}$ satisfying the condition in \eqref{min_unconstrained} follows from the fact that the problem in \eqref{min_unconstrained} is an \textit{unconstrained} discounted MDP problem \cite{puterman2014markov}. The policy $\widetilde{\pi}$ satisfies
\begin{align}\label{opt_lower_bound}
    J(\widetilde{\pi}) \leq \inf_{\pi\in\Xi(\mathcal{M},B)}J(\pi)
\end{align}
since $\Xi(\mathcal{M},B)$$\subseteq$$\Pi(\mathcal{M}')$. Therefore, if $\widetilde{\pi}$$\in$$\Xi(\mathcal{M},B)$, then 
$\widetilde{\pi}$ is an optimal policy for the problem in \eqref{main_objective}. However, in general, we have $\widetilde{\pi}$$\not\in$$\Xi(\mathcal{M},B)$, in which case we need to perturb the policy $\widetilde{\pi}$ in a certain way to obtain an $\epsilon$-optimal solution to the problem in \eqref{main_objective}.

\subsection{Perturbations to maximize reachability}

Given the stationary policy $\widetilde{\pi}$, for each state $s$$\in$$S$, let $\mathcal{A}_{act}(s)$ be the set of actions that are taken by a non-zero probability, i.e.,
\begin{align*}
    \mathcal{A}_{act}(s) \coloneqq \{a\in \mathcal{A}_{\max}(s) \ |\ \widetilde{\pi}(s,a)>0\}.
\end{align*}
Similarly, let $ \mathcal{A}_{pass}(s)$$\coloneqq$$\mathcal{A}_{\max}(s)\backslash\mathcal{A}_{act}(s)$.
We define the perturbed stationary policy $\widetilde{\pi}'$$\in$$\Pi^S(\mathcal{M}')$ as
\begin{align}\label{perturbed_policy_def}
    \widetilde{\pi}'(s,a)\coloneqq\begin{cases} \epsilon' & \text{if} \ \ a\in \mathcal{A}_{pass}(s)\\
    \widetilde{\pi}(s,a)-\epsilon'\frac{\lvert \mathcal{A}_{pass}(s)\rvert}{\lvert \mathcal{A}_{act}(s) \rvert} &  \text{otherwise},
    \end{cases}
\end{align}
where $\epsilon'$$>$$0$ is a sufficiently small constant such that $\widetilde{\pi}'(s,a)$$\geq$$0$ for all $s$$\in$$S$ and $a$$\in$$\mathcal{A}$. Note that the policy $\widetilde{\pi}'$ is well-defined since $\sum_{a\in\mathcal{A}_{\max}(s)}\widetilde{\pi}'(s,a)$$=$$1$ for all $s$$\in$$S$.


{\setlength{\parindent}{0cm}
\noindent \begin{definition} For an MDP $\mathcal{M}$ and a stationary policy $\pi$$\in$$\Pi^S(\mathcal{M})$, an \textit{induced Markov chain} (MC) $\mathcal{M}_{\pi}$$=$$(S,\mathcal{P}^{\pi})$ is a tuple where the transition function $\mathcal{P}^{\pi}$$:$$S$$\times$$S$$\rightarrow$$[0,1]$ is such that, for all $s,s'$$\in$$S$,
\begin{align*}
    \mathcal{P}^{\pi}_{s,s'}\coloneqq\sum_{a\in\mathcal{A}(s)}\pi(s,a)\mathcal{P}_{s,a,s'}.
\end{align*}
\end{definition}}
\noindent We now prove that $\widetilde{\pi}'$ maximizes the probability of reaching  $B$.
{\setlength{\parindent}{0cm}
\noindent \begin{lemma}\label{feasibility_lemma}
It holds that
$\widetilde{\pi}'$$\in$$\Xi(\mathcal{M},B)$.
\end{lemma}}
\noindent\textbf{Proof:} The policy $\widetilde{\pi}'$ is such that, for all $s$$\in$$S^{\widetilde{\pi}'}_{\rightarrow}$, all actions $a$$\in$$\mathcal{A}_{\max}(s)$ are taken with a non-zero probability. Therefore, in the induced MC $\mathcal{M}'_{\widetilde{\pi}'}$, there exists a path to $B$ from each state $s$$\in$$S^{\widetilde{\pi}'}_{\rightarrow}$$\cap$$S_r$. Then, it follows from Lemma 1 in \cite{teichteil2012stochastic} that there exists a constant $N$$\in$$\mathbb{N}$ such that, for all $M$$\geq$$N$, $\varrho^{\widetilde{\pi}'}[M]$$\not\in$$S^{\widetilde{\pi}'}_{\rightarrow}$$\cap$$S_r$. In other words, the agent eventually leaves the states $s$$\in$$S^{\widetilde{\pi}'}_{\rightarrow}$$\cap$$S_r$ with probability 1. By the construction of the cleaned-up MDP $\mathcal{M}'$, the agent reaches the set $B$ with maximum probability under any policy leaving the set $S^{\widetilde{\pi}'}_{\rightarrow}$$\cap$$S_r$ with probability 1. Hence, we have $\widetilde{\pi}'$$\in$$\Xi(\mathcal{M}',B)$. Then, the result follows as $\Xi(\mathcal{M}',B)$$\subseteq$$\Xi(\mathcal{M},B)$.
$\Box$

We now show that, for an appropriately 
chosen constant $\epsilon'$, the stationary policy $\widetilde{\pi}'$ constitutes an $\epsilon$-optimal policy to the problem in \eqref{main_objective}. With an abuse of notation, for a given $\pi$$\in$$\Pi^S(\mathcal{M}')$, let $\mathcal{P}^{\pi}$$\in$$\mathbb{R}^{\lvert S\rvert\times\lvert S\rvert}$ be the transition matrix of the induced MC $\mathcal{M}_{\pi}$.
Let the vector ${\boldsymbol{\alpha}}$$=$$(\alpha_s)_{s\in S}$$\in$$\mathbb{R}^{\lvert S\rvert}$ be the initial state distribution, i.e., $\alpha_{s}$$\coloneqq$$1$ if $s$$=$$s_1$ and $\alpha_{s}$$\coloneqq$$0$ otherwise. Finally, for $\pi$$\in$$\Pi^S(\mathcal{M}')$, let ${\boldsymbol c}^{\pi}$$\coloneqq$$ (c^{\pi}_s)_{s\in S}$$\in$$\mathbb{R}^{\lvert S\rvert}$ be a vector such that
\begin{align*}
    c^{\pi}_s\coloneqq\sum_{a\in\mathcal{A}_{\max}(s)}\pi(s,a)c(s,a).
\end{align*}
It is known \cite{altman1999constrained,puterman2014markov} that, for any $\pi$$\in$$\Pi^S(\mathcal{M}')$, we have
\begin{align}\label{total_cost_matrix}
    J(\pi) = {\boldsymbol{\alpha}}^T(\mathbb{I}-\beta\mathcal{P}^{\pi} )^{-1}{\boldsymbol c}^{\pi}
\end{align}
where $(\cdot)^T$ denotes the transpose operation, and $\mathbb{I}$$\in$$\mathbb{R}^{\lvert S\rvert \times\lvert  S\rvert}$ is the identity matrix. Moreover, it follows from \eqref{perturbed_policy_def} that
\begin{align*}
    \mathcal{P}^{\widetilde{\pi}'} = \mathcal{P}^{\widetilde{\pi}}+\epsilon' M \ \ \ \text{and} \ \ \
    \boldsymbol{c}^{\widetilde{\pi}'} =  \boldsymbol{c}^{\widetilde{\pi}}+\epsilon' {\boldsymbol{v}}
\end{align*}
where $M$$\in$$\mathbb{R}^{\lvert S\rvert \times\lvert S\rvert}$ and ${\boldsymbol v}$$\coloneqq$$(v_s)_{s\in S}$$\in$$\mathbb{R}^{\lvert S\rvert}$ are, respectively, the perturbation matrix and the perturbation vector satisfying
\begin{align*}
    M(s,s')\coloneqq\sum_{a\in \mathcal{A}_{pass}(s)}\mathcal{P}_{s,a,s'}-\frac{\lvert \mathcal{A}_{pass}(s) \rvert}{\lvert \mathcal{A}_{act}(s) \rvert}\sum_{a\in \mathcal{A}_{act}(s)}\mathcal{P}_{s,a,s'},\\
    v_s \coloneqq \sum_{a\in \mathcal{A}_{pass}(s)}c(s,a)-\frac{\lvert \mathcal{A}_{pass}(s) \rvert}{\lvert \mathcal{A}_{act}(s) \rvert}\sum_{a\in \mathcal{A}_{act}(s)}c(s,a).
\end{align*}

{\setlength{\parindent}{0cm}
\noindent \begin{thm}
For any given $\epsilon$$>$$0$, the policy $\widetilde{\pi}'$ defined in \eqref{perturbed_policy_def} is an $\epsilon$-optimal policy for the problem in \eqref{main_objective} if $\epsilon'$$>$$0$ is chosen such that
\begin{align}\label{epsilon_upper_bound}
    \epsilon' \leq \frac{\epsilon}{\gamma_1 + \gamma_2}
\end{align}
where
$\gamma_1\coloneqq\beta\boldsymbol{\alpha}^T(\mathbb{I}-\beta\mathcal{P}^{\widetilde{\pi}} )^{-1}M(\mathbb{I}-\beta\mathcal{P}^{\widetilde{\pi}'} )^{-1}{\boldsymbol c}^{\widetilde{\pi}}$ and $\gamma_2\coloneqq{\boldsymbol{\alpha}}^T(\mathbb{I}-\beta\mathcal{P}^{\widetilde{\pi}'} )^{-1}{\boldsymbol v}$.
\end{thm}}
\noindent \textbf{Proof:} We show in Lemma \ref{feasibility_lemma} that the policy $\widetilde{\pi}'$ is feasible for the problem in \eqref{main_objective}. Here, we show that the cost-to-go $J(\widetilde{\pi}')$ is at most $\epsilon$ larger than the minimum achievable one.

Using \eqref{total_cost_matrix} and the definition of $\gamma_2$, we have
\begin{subequations}
\begin{align}
    J(\widetilde{\pi}')&= {\boldsymbol{\alpha}}^T(\mathbb{I}-\beta\mathcal{P}^{\widetilde{\pi}'} )^{-1}{\boldsymbol c}^{\widetilde{\pi}'}\\
    &= {\boldsymbol{\alpha}}^T(\mathbb{I}-\beta\mathcal{P}^{\widetilde{\pi}'} )^{-1}({\boldsymbol c}^{\widetilde{\pi}}+\epsilon' {\boldsymbol v})\\ \label{last_cost_eq}
    &= {\boldsymbol{\alpha}}^T(\mathbb{I}-\beta\mathcal{P}^{\widetilde{\pi}'} )^{-1}{\boldsymbol c}^{\widetilde{\pi}}+ \epsilon' \gamma_2.
\end{align}
\end{subequations}
It follows from  equation (26) in \cite{henderson} that
\begin{subequations}
\begin{align}
   (\mathbb{I}-\beta\mathcal{P}^{\widetilde{\pi}'} )^{-1}&= (\mathbb{I}-\beta\mathcal{P}^{\widetilde{\pi}}-\epsilon' \beta M )^{-1}\\
   &=(\mathbb{I}-\beta\mathcal{P}^{\widetilde{\pi}})^{-1}\nonumber\\ \label{matrix_inversion}
   & \ \ \  +\epsilon'\beta(\mathbb{I}-\beta\mathcal{P}^{\widetilde{\pi}})^{-1} M(\mathbb{I}-\beta\mathcal{P}^{\widetilde{\pi}'})^{-1}.
\end{align}
\end{subequations}
Plugging \eqref{matrix_inversion} into \eqref{last_cost_eq}, we obtain $J(\widetilde{\pi}')$$=$$J(\widetilde{\pi})$$+$$\epsilon' (\gamma_1$$+$$\gamma_2)$.
Using \eqref{epsilon_upper_bound} and the definition of $\widetilde{\pi}$ given in \eqref{min_unconstrained}, we conclude that
\begin{align}\label{last_ineq_epsilon}
    J(\widetilde{\pi}') \leq \min_{\pi\in\Pi(\mathcal{M}')}J(\pi)+\epsilon.
\end{align}
Then, the result follows as $\Xi(\mathcal{M},B)$$\subseteq$$\Pi(\mathcal{M}')$. $\Box$

We conclude this section by summarizing the three main steps of the efficient synthesis of an $\epsilon$-optimal stationary policy for the problem in \eqref{main_objective}. First, obtain the cleaned-up MDP $\mathcal{M}'$ by removing all actions $a$$\in$$\mathcal{A}(s)\backslash\mathcal{A}_{\max}(s)$ where $\mathcal{A}_{\max}(s)$ is as defined in \eqref{A_max_def}. Second, synthesize a stationary deterministic policy $\widetilde{\pi}$$\in$$\Pi^{SD}(\mathcal{M}')$ that satisfies \eqref{min_unconstrained} via linear programming. Finally, obtain an $\epsilon$-optimal policy $\widetilde{\pi}'$ given in \eqref{perturbed_policy_def} by choosing $\epsilon'$ as shown in \eqref{epsilon_upper_bound}.

\section{A Necessary and Sufficient Condition for the Existence of Optimal Policies}
In the previous sections, we showed that an optimal policy solving the problem in \eqref{main_objective} may not exist, but an $\epsilon$-optimal stationary policy can be synthesized efficiently. To complete the analysis, we now provide an efficiently verifiable necessary and sufficient condition for the existence of an optimal policy that solves the problem in \eqref{main_objective}.

The following result is an immediate consequence of the inequality in \eqref{last_ineq_epsilon}, which shows that the optimal value of the problem in \eqref{main_objective} coincides with the optimal value of the problem in \eqref{min_unconstrained}.
{\setlength{\parindent}{0cm}
\noindent \begin{Corollary}\label{corol_1}
The following equality holds:
\begin{align}\label{corollary_equality}
    \inf_{\pi\in\Xi(\mathcal{M},B)}J(\pi) = \min_{\pi\in\Pi(\mathcal{M}')}J(\pi).
\end{align}
\end{Corollary}}
{\setlength{\parindent}{0cm}
\noindent \begin{thm}\label{thm_2}
There exists an optimal policy $\pi^{\star}$$\in$$\Pi(\mathcal{M})$ that satisfies the condition in \eqref{main_objective} if and only if there exists a policy $\widetilde{\pi}$$\in$$\Xi(\mathcal{M}',B)$ that satisfies the condition in \eqref{min_unconstrained}.
\end{thm}}
\noindent\textbf{Proof:} ($\Leftarrow$) Suppose that there exists a policy $\widetilde{\pi}$$\in$$\Xi(\mathcal{M}',B)$ that satisfies the condition in \eqref{min_unconstrained}. Then, $\widetilde{\pi}$ is an optimal policy satisfying the condition in \eqref{main_objective} due to \eqref{opt_lower_bound}.

($\Rightarrow$) Suppose that there exists an optimal policy $\pi^{\star}$$\in$$\Pi(\mathcal{M})$ that satisfies the condition in \eqref{main_objective}. Since $\pi^{\star}$ is a feasible solution to the problem in \eqref{main_objective}, it satisfies $\pi^{\star}$$\in$$\Xi(\mathcal{M}',B)$. Additionally, it follows from \eqref{corollary_equality} that $\pi^{\star}$ satisfies the condition in \eqref{min_unconstrained}. $\Box$

Theorem \ref{thm_2} establishes that, for a given MDP, the existence of an optimal policy solving the problem in \eqref{main_objective} can be verified in three steps as follows. First, construct the cleaned-up MDP $\mathcal{M}'$ as described in Section \ref{cleaned_up_sec}. Second, find all policies that minimize the total discounted cost on the cleaned-up MDP $\mathcal{M}'$, i.e., policies that satisfy the condition in \eqref{min_unconstrained}. Finally, check whether any of these policies belong to the set $\Xi(\mathcal{M}',B)$. In Section \ref{cleaned_up_sec}, we show that the first step can be performed efficiently. In what follows, we show that the second and the third steps of the above procedure can also be performed efficiently.

Let ${\bf{y}}$$:=$$(y_s)_{s\in S}$$\in$$\mathbb{R}^{\lvert S\rvert}$ be a vector such that
\begin{align*}
    y_s:= \min\limits_{{a\in \mathcal{A}_{\max}(s)}}\Bigg\{ c(s,a)+\beta \sum_{s'\in S}\mathcal{P}_{s,a,s'} y_{s'} \Bigg\}.
\end{align*}
It is known \cite{puterman2014markov} that each element $y_s$ of ${\bf{y}}$ corresponds to the minimum total discounted cost accumulated along the paths that start from the state $s$$\in$$S$. The vector ${\bf{y}}$ can be efficiently computed via value iteration or linear programming \cite{puterman2014markov}. The following result provides a necessary and sufficient condition for a policy to minimize the total discounted cost in an MDP in the absence of constraints.
{\setlength{\parindent}{0cm}
\noindent \begin{prop} \label{propositon_cost} \cite[Chapter 6]{puterman2014markov}
A policy $\widetilde{\pi}$$=$$(\widetilde{d}_1,\widetilde{d}_2,\widetilde{d}_3,\ldots)$$\in$$\Pi(\mathcal{M}')$ satisfies the condition in \eqref{min_unconstrained} if and only if, for all $s$$\in$$ S_{\rightarrow}^{\widetilde{\pi}}$,
\begin{align}\label{necc_opt}
    y_s= c(s,a)+\beta \sum_{s'\in S}\mathcal{P}_{s,a,s'} y_{s'}
\end{align}
for all $a$$\in$$\mathcal{A}_{\max}(s)$ satisfying $\widetilde{d}_t(s,a)$$>$$0$ for some $t$$\in$$\mathbb{N}$.
\end{prop}}

Given a cleaned-up MDP $\mathcal{M}'$, we obtain the \textit{modified MDP} $\overline{\mathcal{M}}'$ by removing all actions $a$$\in$$\mathcal{A}_{\max}(s)\backslash \mathcal{A}_{opt}(s)$ from each $s$$\in$$S$, where
\begin{align*}
    \mathcal{A}_{opt}(s):=\Bigg\{a\in\mathcal{A}_{\max}(s)\  \Big| \ y_s= c(s,a)+\beta \sum_{s'\in S}\mathcal{P}_{s,a,s'} y_{s'}\Bigg\}.
\end{align*}
Then, it follows from Proposition \ref{propositon_cost} that the policies on the modified MDP $\overline{\mathcal{M}}'$, and only them, minimize the total discounted cost on $\mathcal{M}'$.

Finally, we verify whether there exists a policy $\pi$$\in$$\Pi(\overline{\mathcal{M}}')$ on the modified MDP such that $\pi$$\in$$\Xi(\mathcal{M}',B)$. Let ${\overline{\bf{x}}}$$:=$$(\overline{x}_s)_{s\in S}$$\in$$\mathbb{R}^{S}$ be a vector such that the maximum in \eqref{x_vector} is taken over the set $\mathcal{A}_{opt}(s)$ instead of the set $\mathcal{A}(s)$. Each element of ${\overline{\bf{x}}}$ corresponds to the probability of reaching the target set $B$ from the state $s$$\in$$S$ in $\overline{\mathcal{M}}'$, which can be computed efficiently via linear programming.
{\setlength{\parindent}{0cm}
\noindent \begin{prop}  \label{prop_3_reach}
There exists a policy $\pi$$\in$$\Pi(\overline{\mathcal{M}}')$ such that $\pi$$\in$$\Xi(\mathcal{M}',B)$ if and only if $\overline{x}_{s_1}$$=$$x_{s_1}$.
\end{prop}}
\noindent \textbf{Proof:}  ($\Rightarrow$) If $\pi$$\in$$\Pi(\overline{\mathcal{M}}')$ satisfies the condition $\pi$$\in$$\Xi(\mathcal{M}',B)$, then the result follows from Proposition \ref{prop_1_statement} and the fact that $\overline{x}_{s_1}$$\leq$$x_{s_1}$.

($\Leftarrow$) If $\overline{x}_{s_1}$$=$$x_{s_1}$, then a policy $\pi$ satisfying $\pi$$\in$$\Xi(\mathcal{M}',B)$ can be constructed as explained in Lemma 10.102 of \cite{Model_checking}. Specifically, on the modified MDP, we first remove the actions from each state $s$$\in$$S_r$ that do not satisfy the condition $\overline{x}_s$$=$$\sum_{s'\in S}\mathcal{P}_{s,a,s'}\overline{x}_{s'}$. On the graph corresponding to the resulting MDP, let $T_{\min}(s)$ be the length of the shortest path from the state $s$$\in$$S_r$ to the target set $B$ (see Section \ref{app_algo_sec} for details). We construct a policy $\pi$$\in$$\Xi(\mathcal{M}',B)$ by choosing an action in each state $s$$\in$$S_r$ under which the agent transitions with a nonzero probability to a state $s'$ such that $T_{\min}(s')$$<$$T_{\min}(s)$.  $\Box$

We conclude this section by noting that, once the existence of an optimal policy solving the problem in \eqref{main_objective} is verified using the procedure explained above, one can synthesize such a policy through the procedure described in the proof of Proposition \ref{prop_3_reach}.

\section{An Analysis Over Deterministic Policies}\label{det_analysis_sec}
In many applications, it is desirable to optimize performance using deterministic policies due to their low computational requirements and ease of implementation in distributed systems \cite{paruchuri2004towards}. Accordingly, in this section, we focus on stationary deterministic policies and consider the problem of synthesizing a policy $\pi_D^{\star}$ such that
\begin{align}\label{det_objective_2}
  \pi_D^{\star} \in \arg \min_{\pi\in\Xi^{SD}(\mathcal{M},B)}\  J(\pi)
\end{align}
where the set $\Xi^{SD}(\mathcal{M},B)$$\subseteq$$\Pi^{SD}(\mathcal{M})$ denotes the set of stationary deterministic policies that maximize the probability of reaching the target set $B$. Specifically, it holds $\pi$$\in$$\Xi^{SD}(\mathcal{M},B)$ if and only if $\pi$$\in$$\Pi^{SD}(\mathcal{M})$ satisfies the condition in \eqref{max_reac_def}. Since the set $\Xi^{SD}(\mathcal{M},B)$ is finite, an optimal policy for the problem in \eqref{det_objective_2} always exists.

\subsection{A complexity result} We now show that the synthesis of a policy that satisfies the condition in \eqref{det_objective_2} is, in general, intractable. We remark that the synthesis of a stationary deterministic policy that minimizes the total discounted cost in MDPs subject to total discounted reward constraints is known to be NP-hard \cite{feinberg2000constrained}. The problem in \eqref{det_objective_2} involves a probabilistic constraint which, in general, cannot be represented as a discounted criterion without changing the feasible policy space. Therefore, 
none of the existing complexity results on constrained MDPs apply to 
the problem in \eqref{det_objective_2}.

We prove the result by a reduction from the Hamiltonian path problem (HAMPATH) which is known to be NP-complete \cite{sipser1996introduction}. Let $G$$=$$(V,E)$ be a directed graph (digraph) where $V$ is a finite set of vertices and $E$$\subseteq$$V$$\times$$V$ is a finite set of edges. For a given digraph $G$, a \textit{finite path} $v_1v_2\ldots v_n$ of length $n$$\in$$\mathbb{N}$ from vertex $v_1$ to $v_n$ is a sequence of vertices such that $(v_k,v_{k+1})$$\in$$E$ for all $1$$\leq$$k$$<$$n$.
{\setlength{\parindent}{0cm}
\noindent \begin{definition} (HAMPATH)
Given a digraph $G$$=$$(V,E)$ and an origin-destination pair $(o,d)$$\in$$V$$\times$$V$, decide whether there exists a Hamiltonian path on $G$, i.e., a finite path from the vertex $o$ to the vertex $d$ that visits each vertex $v$$\in$$V$ exactly once.
\end{definition} \begin{thm}
For $K$$\in$$\mathbb{N}$, deciding whether there exists a policy $\pi$$\in$$\Xi^{SD}(\mathcal{M},B)$ such that $J(\pi)$$\leq$$K$ is NP-complete.
\end{thm}}
\noindent \textbf{Proof:} The decision problem is in NP since for any given policy $\pi$$\in$$\Xi^{SD}(\mathcal{M},B)$, we can verify whether $J(\pi)$$\leq$$K$ in polynomial-time using the formula in \eqref{total_cost_matrix}.

To show the NP-hardness, we reduce an arbitrary HAMPATH instance to an instance of the problem in \eqref{det_objective_2}.


Given a graph $G$ and a pair $(o,d)$, we construct an instance of the problem in \eqref{det_objective_2} such that the agent starts from $o$ and aims to reach $d$. We define the MDP $\mathcal{M}$ as follows. The set of states are $S$$=$$V$, the initial state is $s_1$$=$$o$, and the target set is $B$$=$$\{d\}$. For a state $v$$\in$$V$, we associate an action $a_{v'}$ for each edge $(v,v')\in E$, i.e., $\mathcal{A}(v)$$=$$\{a_{v'}$$\in$$\mathcal{A} : (v,v')$$\in$$E \}$. The transition function $\mathcal{P}$ is such that $\mathcal{P}_{v,a_{v'},v'}$$=$$1$ if $v$$\in$$V\backslash B$, and $\mathcal{P}_{v,a_{v'},v}$$=$$1$ if $v$$\in$$B$. Finally, we define the cost function $c$ such that
\begin{align*}
    c(v,a_{v'})=\begin{cases} K/\beta^{\lvert V \rvert -1} &\ \  \text{if} \ \ v'=d\\
    0 &\ \  \text{otherwise}.
    \end{cases}
\end{align*}

We now show that, there exists a policy $\pi$$\in$$\Xi^{SD}(\mathcal{M},B)$ such that $J(\pi)$$\leq$$K$ if and only if there exists a Hamiltonian path on the graph $G$ with the origin-destination pair $(o,d)$.

Suppose that there is a Hamiltonian path $v_1v_2\ldots v_{\lvert V\rvert}$ where $v_1$$=$$o$ and $v_{\lvert V\rvert}$$=$$d$. Consider the policy $\pi$$\in$$\Pi^{SD}(\mathcal{M})$ such that $\pi(v_k,a_{v_{k+1}})$$=$$1$ for all $k$$\in$$\mathbb{N}$ such that $k$$<$$\lvert V\rvert$. Under the policy $\pi$, the agent reaches the target set $B$ with probability one; hence, it holds $\pi$$\in$$\Xi^{SD}(\mathcal{M},B)$. Moreover, the agent reaches the set $B$ in exactly $\lvert V\rvert$$-$$1$ steps; hence, it holds that $J(\pi)$$=$$K$.

Suppose that there exists a policy $\pi$$\in$$\Xi^{SD}(\mathcal{M},B)$ such that $J(\pi)$$\leq$$K$. Then, under the policy $\pi$, the agent reaches the target set $B$ with probability one in at least $\lvert V \rvert$$-$$1$ steps. Moreover, since $\pi$ is a deterministic policy, the agent visits each state at most once because, otherwise, the probability of reaching the target set $B$ must be zero. The combination of the above arguments imply that the agent reaches the target set $B$ in $\lvert V \rvert$$-$$1$ steps by visiting each state $v$$\in$$V\backslash B$ exactly once. By definition, such a path constitutes a Hamiltonian path. $\Box$

\subsection{An exact algorithm}\label{exact_algorithm_sec}
In this section, we solve the problem in \eqref{det_objective_2} by formulating it as a mixed-integer linear program (MILP). The MILP formulation is an extension of the results in \cite{dolgov2005stationary,kalagarla2020synthesis}. In \cite{dolgov2005stationary}, the authors synthesize stationary deterministic policies for MDPs involving multiple discount factors each of which is \textit{strictly} less than one. In \cite{kalagarla2020synthesis}, an MILP formulation is presented for the synthesis of stationary deterministic policies under which an agent satisfies a temporal logic specification \textit{with probability one} while minimizing its total discounted cost. Here, we formulate an MILP for the synthesis of a policy that reaches a desired set of target states with maximum probability (which is potentially less than one) while minimizing the total discounted cost.

Consider the following MILP:
\begin{subequations}
\begin{flalign}\label{MILP_start}
    &\underset{\substack{\lambda^1_{s,a}, \lambda^2_{s,a}, \Delta_{s,a}}}{\text{minimize}} \sum_{s\in S}\sum_{a\in \mathcal{A}} c(s,a) \lambda^1_{s,a}\\
    &\text{subject to:} \nonumber \\ \label{MILP_cons_1}
    &\ \ \ \sum_{a\in \mathcal{A}}\lambda^1_{s,a}-\beta\sum_{s'\in S}\sum_{a\in\mathcal{A}} \mathcal{P}_{s',a,s}\lambda^1_{s',a}= \alpha_s, \  \forall s\in S\\ \label{MILP_cons_2}
    &\ \ \ \sum_{a\in \mathcal{A}}\lambda^2_{s,a}-\sum_{s'\in S}\sum_{a\in\mathcal{A}} \mathcal{P}_{s',a,s}\lambda^2_{s',a}= \alpha_s, \ \ \  \forall s\in S_r\\ \label{MILP_cons_3}
    &\ \ \ \sum_{s\in S_r}\sum_{a\in \mathcal{A}}\lambda^2_{s,a}r(s,a) = x_{s_1}\\ \label{MILP_cons_4}
    & \ \ \ \lambda^k_{s,a}/M \leq \Delta_{s,a}, \quad \forall k \in \{1,2\}, \   \forall s \in S, \ \forall a\in \mathcal{A}\\ \label{MILP_cons_5}
    & \ \ \ \sum_{a\in \mathcal{A}}\Delta_{s,a} \leq 1, \qquad \qquad \qquad \qquad \qquad \quad \forall s\in S\\
    &\ \ \ \lambda^1_{s,a}\geq 0, \lambda^2_{s,a}\geq 0, \Delta_{s,a}\in \{0,1\} \label{MILP_end}
\end{flalign}
\end{subequations}
where $M$ is a large constant whose precise value will be discussed shortly. The vector ${\boldsymbol{\alpha}}$$=$$(\alpha_s)_{s\in S}$ denotes the initial state distribution, i.e., $\alpha_s$$=$$1$ if $s$$=$$s_1$ and $\alpha_s$$=$$0$ otherwise; $x_{s_1}$ denotes the maximum probability of reaching the target set $B$ from the initial state and can be computed via linear programming as discussed in Section \ref{cleaned_up_sec}. Finally, $r$$:$$S$$\times$$\mathcal{A}$$\rightarrow$$\mathbb{R}_{\geq 0}$ is a function such that
\begin{align*}
    r(s,a):=\begin{cases} \sum_{s'\in B}\mathcal{P}_{s,a,s'} & \ \text{if}\ \  s \in S_r \\
    0 & \ \text{otherwise}.
    \end{cases}
\end{align*}

There are $(\lvert S\rvert \lvert \mathcal{A}\rvert)^2$ continuous and $\lvert S\rvert \lvert \mathcal{A}\rvert$ integer variables in \eqref{MILP_start}-\eqref{MILP_end}. The sets $\{\lambda^1_{s,a} $$\geq$$0 : s$$\in$$S, a$$\in$$\mathcal{A}\}$ and $\{\lambda^2_{s,a}$$\geq$$0 : s$$\in$$S, a$$\in$$\mathcal{A}\}$ of continuous variables correspond, respectively, to the \textit{discounted} and \textit{undiscounted} occupation measures \cite{altman1999constrained,etessami2007multi}. The set $\{\Delta_{s,a}$$\in$$\{0,1\} : s$$\in$$S, a$$\in$$\mathcal{A}\}$ of integer variables correspond to the deterministic actions taken by the agent.

The constraints in \eqref{MILP_cons_1} and \eqref{MILP_cons_2} represent the balance equations for the discounted and undiscounted occupation measures, respectively \cite{altman1999constrained}. The constraint in \eqref{MILP_cons_3} ensures that the agent reaches the target set $B$ with maximum probability $x_{s_1}$. Finally, the constraints in \eqref{MILP_cons_4}-\eqref{MILP_cons_5} ensure that the agent follows a deterministic policy.

Let $f$$:$$S$$\times$$\mathcal{A}$$\rightarrow$$\mathbb{R}_{\geq 0}$ be a function such that $f(s,a)$$:=$$1$ for all $s$$\in$$S_r$ and $f(s,a)$$:=$$0$ otherwise. Moreover, let $M^{\star}$ be a constant such that
\begin{align*}
    M^{\star} := \max_{\pi\in \Xi^{SD}(\mathcal{M},B)} \mathbb{E}^{\pi}\Bigg[\sum_{t=1}^{\infty} f(s_t,a_t)\Bigg].
\end{align*}
$M^{\star}$ corresponds to the maximum expected number of steps taken under any stationary deterministic policy before reaching the set $B$ with maximum probability. $M^{\star}$ is finite since all policies $\pi$$\in$$ \Xi^{SD}(\mathcal{M},B)$ induce MCs in which the set $S_r$$\cap$$S^{\pi}_{\rightarrow}$ consists only of transient states.{\setlength{\parindent}{0cm}
\noindent \begin{thm}\label{MILP_thm}
Let the constant $M$ in \eqref{MILP_cons_4} is chosen such that $M$$\geq$$M^{\star}$. Then, for $K$$\in$$\mathbb{R}$, there exists a policy $\pi$$\in$$\Xi^{SD}(\mathcal{M},B)$ such that $J(\pi)$$\leq$$K$ if and only if the optimal value of the problem in \eqref{MILP_start}-\eqref{MILP_end} is less than or equal to $K$. Furthermore, an optimal policy $\pi_D^{\star}$$\in$$\Xi^{SD}(\mathcal{M},B)$ such that $J(\pi_D^{\star})$$\leq$$K$ can be obtained from the optimal variables $\{\Delta_{s,a}^{\star} : s$$\in$$S,  a$$\in$$\mathcal{A}\}$ of the problem in \eqref{MILP_start}-\eqref{MILP_end} by the following rule:
\begin{align}\label{opt_policy_construct}
    \pi_D^{\star}(s,a)=\begin{cases} \Delta_{s,a}^{\star}, & \text{if} \ \ \sum_{a\in \mathcal{A}} \Delta_{s,a}^{\star}=1\\
    \text{arbitrary}, & \text{otherwise.}
    \end{cases}
\end{align}
\end{thm}}
\noindent \textbf{Proof:} ($\Leftarrow$) Suppose that there exists a stationary deterministic policy $\pi$$\in$$\Xi^{SD}(\mathcal{M},B)$ such that $J(\pi)$$\leq$$K$. We set $\Delta_{s,a}$$=$$\pi(s,a)$,
\begin{align*}
    \lambda^1_{s,a}&=\sum_{t=1}^{\infty}\beta^{t-1}\text{Pr}^{\pi}\Big(S_t = s, A_t = a| S_1 = s_1\Big), \ \text{and}\\
     \lambda^2_{s,a}&=\begin{cases} \sum_{t=1}^{\infty}\text{Pr}^{\pi}\Big(S_t = s, A_t = a | S_1= s_1\Big), & \text{if} \ \ s\in S_r\\
     0, & \text{otherwise}.
     \end{cases}
\end{align*}
 $\text{Pr}^{\pi}(S_t$$=$$s, A_t$$=$$a | S_1$$=$$s_1)$ denotes the probability with which the state-action pair $(s,a)$ is occupied in the induced MC $\mathcal{M}^{\pi}$. The variables $\lambda^1_{s,a},\lambda^2_{s,a}$ are finite for all $s$$\in$$S$ and $a$$\in$$\mathcal{A}$ since the states $s$$\in$$S_r$$\cap$$S^{\pi}_{\rightarrow}$ are transient in the induced MC $\mathcal{M}^{\pi}$.

It can be shown that the above choices of the variables satisfy the balance equations in \eqref{MILP_cons_1}-\eqref{MILP_cons_2} (see, e.g., Theorem 3.1 in \cite{altman1999constrained} for a similar derivation). The variables $\lambda^2_{s,a}$ also satisfy the constraint in \eqref{MILP_cons_3}, which follows from the fact (see Theorem 10.15 in \cite{Model_checking}) that, for any $\pi$$\in$$\Pi(\mathcal{M})$, we have
\begin{align}\label{equivalence_reach}
     \text{Pr}_{\mathcal{M}}^{\pi}(Reach[B]) = \mathbb{E}^{\pi}\Bigg[\sum_{t=1}^{\infty} r(s_t,a_t)\Bigg].
\end{align}
The satisfaction of the constraint in \eqref{MILP_cons_4} follows from the definition of $M^{\star}$. The constraints in \eqref{MILP_cons_5} and \eqref{MILP_end} are satisfied by construction. Finally, it follows from  classical MDP theory, (see, e.g., \cite[Chapter 6]{puterman2014markov}), that these variables attain an objective value that is equal to $J(\pi)$. Since we constructed a feasible solution with the objective value $J(\pi)$, we conclude that the optimal value of the problem in \eqref{MILP_start}-\eqref{MILP_end} is less than or equal to $K$.

($\Rightarrow$) Suppose that the optimal value of the problem in \eqref{MILP_start}-\eqref{MILP_end} is less than or equal to $K$. From the optimal variables, we construct a policy $\pi_D^{\star}$ through the formula in \eqref{opt_policy_construct}. It follows from Proposition 2 in \cite{dolgov2005stationary} that $\pi_D^{\star}$$\in$$\Pi^{SD}(\mathcal{M})$ and $J(\pi_D^{\star})$$\leq$$K$. Due to the equality in \eqref{equivalence_reach}, we also have $\pi^{\star}$$\in$$\Xi^{SD}(\mathcal{M},D)$, which concludes the proof. $\Box$

Theorem \ref{MILP_thm} establishes that, if we choose $M$$\geq$$M^{\star}$, we can obtain a solution to the problem in \eqref{det_objective_2} by solving the MILP in \eqref{MILP_start}-\eqref{MILP_end} and constructing a policy through the formula in \eqref{opt_policy_construct}. It can be shown by a reduction from the Hamiltonian path problem that the computation of the constant $M^{\star}$ is NP-hard in general.  For MDPs with deterministic transitions, we have $M^{\star}$$\leq$$\lvert S\rvert$; hence we set $M$$=$$\lvert S\rvert$. Although a finite upper bound on $M^{\star}$ for MDPs with stochastic transitions can be computed efficiently, we omit the details as the procedure is more involved and requires one to analyze the structure of the MDP. For such MDPs, we simply set $M$$=$$k\lvert S \rvert$ for some large $k$$\in$$\mathbb{N}$.


\subsection{An approximation algorithm}\label{app_algo_sec}
We now present an algorithm to obtain an \textit{approximate} solution to the NP-hard policy synthesis problem given in \eqref{det_objective_2}. We make the following assumption throughout this section.

\noindent \textbf{Assumption:} The cost function $c$$:$$S$$\times$$\mathcal{A}$$\rightarrow$$\mathbb{R}_{\geq 0}$ satisfies $c(s,a)$$=$$0$ for all $s$$\in$$B\cup S_0$.

The above assumption states that the agent incurs no cost at the target states and at the states from which there is no path to the target states. The above assumption typically holds in practice since costs are incurred only until a task is completed or failed to be completed.

The main idea behind the approximation algorithm is that the discounted immediate cost $\beta^{t-1} c(s,a)$ incurred at a state-action pair $(s,a)$ at step $t$$\in$$\mathbb{N}$ is upper bounded by $\beta^{T_{\min}(s)-1} c(s,a)$ where $T_{\min}(s)$ is the minimum number of steps taken to reach the state $s$. Using the derived upper bound, we define a surrogate objective function and formulate a constrained MDP problem for which we synthesize an optimal deterministic policy via linear programming.

For an MDP $\mathcal{M}$, let $G_{\mathcal{M}}$$=$$(S,E_{\mathcal{M}})$ be a digraph where $S$ is the set of vertices and $E_{\mathcal{M}}$ is the set of edges such that
\begin{align}\label{edge_def}
    E_{\mathcal{M}}:=\Bigg\{(s,s')\in S\times S : \sum_{a\in\mathcal{A}(s)}\mathcal{P}_{s,a,s'}>0  \Bigg\}.
\end{align}
On the digraph $G_{\mathcal{M}}$, let $\mathbb{PATH}_n(s)$ be the set of finite paths $v_1v_2\ldots v_n$ of length $n$$\in$$\mathbb{N}$ such that $v_1$$=$$s_1$ and $v_n$$=$$s$, i.e., the set of finite paths that reach the state $s$ starting from the initial state $s_1$. Then, on the MDP $\mathcal{M}$, the agent can reach the state $s$ in minimum $T_{\min}(s)$ steps where
\begin{align}\label{T_min_def}
    T_{\min}(s) := \min \{n \in \mathbb{N} : \mathbb{PATH}_n(s) \neq \emptyset\}.
\end{align}
Note that $T_{\min}(s)$ can be efficiently computed using standard shortest path algorithms, e.g., Dijkstra's algorithm \cite{dijkstra1959note}.

Let $\widetilde{c}$$:$$S$$\times$$\mathcal{A}$$\rightarrow$$\mathbb{R}_{\geq 0}$ be a modified cost function such that
\begin{align}\label{modified_cost_def}
    \widetilde{c}(s,a) := \beta^{T_{\min}(s)-1}c(s,a),
\end{align}
and $\widetilde{J}$$:$$\Pi(\mathcal{M})$$\rightarrow$$\mathbb{R}_{\geq 0}$ be a surrogate function such that
\begin{align*}
    \widetilde{J}(\pi) := \mathbb{E}^{\pi}\Bigg[\sum_{t=1}^{\infty} \widetilde{c}(s_t,a_t) \Bigg].
\end{align*}
As the approximation algorithm, we propose to synthesize a policy $\pi^{\star}_{D,app}$$\in$$\Pi^{SD}(\mathcal{M})$ such that
\begin{align}\label{det_objective_approx}
  \pi^{\star}_{D,app} \in \arg \min_{\pi\in\Xi^{SD}(\mathcal{M},B)}\  \widetilde{J}(\pi).
\end{align}
In what follows, we first present an efficient method to synthesize a stationary deterministic policy $\pi^{\star}_{D,app}$ that satisfies the condition in \eqref{det_objective_approx}. We then derive an upper bound on the suboptimality of the synthesized policy for the original problem given in \eqref{det_objective_2}.

\subsubsection{Policy synthesis} The problem in \eqref{det_objective_approx} is a total \textit{undiscounted} cost minimization problem subject to a constraint on the total \textit{undiscounted} reward given in \eqref{equivalence_reach}. We recently established the existence of optimal stationary deterministic policies for such problems in \cite{savas2020complexity}, where we also present an efficient method to synthesize optimal policies. For completeness, we provide the developed method below.

We solve two linear programs (LPs) to synthesize the policy $\pi^{\star}_{D,app}$. First, we solve the following LP:
\begin{subequations}
\begin{flalign}\label{LP_start}
    &\underset{\substack{\lambda_{s,a}\geq 0}}{\text{minimize}} \sum_{s\in S_r}\sum_{a\in \mathcal{A}} \widetilde{c}(s,a) \lambda_{s,a}\\
    &\text{subject to:} \nonumber \\  \label{LP_cons_1}
    &\ \ \ \sum_{a\in \mathcal{A}}\lambda_{s,a}-\sum_{s'\in S}\sum_{a\in\mathcal{A}} \mathcal{P}_{s',a,s}\lambda_{s',a}= \alpha_s, \  \forall s\in S_r\\
    &\ \ \ \sum_{s\in S_r}\sum_{a\in \mathcal{A}}\lambda_{s,a}r(s,a) = x_{s_1}.\label{LP_end}
\end{flalign}
\end{subequations}
The above LP computes a set $\{\lambda_{s,a}$$\geq$$0: s$$\in$$S, a$$\in$$\mathcal{A}\}$ of occupancy measures from which one can synthesize a (randomized) policy $\pi$ that minimizes the surrogate function $\widetilde{J}(\pi)$ over $\Pi(\mathcal{M})$ while satisfying the maximum reachability constraint. Let $v^{\star}$ be the optimal value of the LP in \eqref{LP_start}-\eqref{LP_end}. Next, we solve the following LP.
\begin{subequations}
\begin{flalign}\label{LP2_start}
    &\underset{\substack{\lambda_{s,a}\geq 0}}{\text{minimize}} \ \ \ \sum_{s\in S}\sum_{a\in \mathcal{A}} \lambda_{s,a}\\
    &\text{subject to:}  \ \ \sum_{s\in S_r}\sum_{a\in \mathcal{A}}\lambda_{s,a}\widetilde{c}(s,a) = v^{\star}\\
    &\qquad \qquad  \ \ \ \eqref{LP_cons_1}-\eqref{LP_end}. \label{LP2_end}
\end{flalign}
\end{subequations}
Let $\{\lambda^{\star}_{s,a}$$\geq$$0: s$$\in$$S, a$$\in$$\mathcal{A}\}$ be the set of optimal variables for the LP in \eqref{LP2_start}-\eqref{LP2_end}. Moreover, for a given state $s$$\in$$S$, let $\mathcal{A}^{\star}(s)$$:=$$\{a$$\in$$\mathcal{A}(s): \lambda^{\star}_{s,a}$$>$$0\}$ be the set of optimal actions. Then, a stationary deterministic policy $\pi^{\star}_{D,app}$$\in$$\Pi^{SD}(\mathcal{M})$ satisfying the condition in \eqref{det_objective_approx} can be synthesized from this set by choosing
\begin{align}\label{opt_policy_rule_approx}
   \pi^{\star}_{D,app}(s,a)=1 \ \text{for an arbitrary} \ a\in\mathcal{A}^{\star}(s).
\end{align}
{\setlength{\parindent}{0cm}
\noindent \begin{prop}\cite{savas2020complexity} A policy $\pi^{\star}_{D,app}$$\in$$\Pi^{SD}(\mathcal{M})$ generated by the rule in \eqref{opt_policy_rule_approx} satisfies the condition in \eqref{det_objective_approx}.
\end{prop}}

The proposed approximation algorithm can be summarized as follows. First, construct the modified cost function given in \eqref{modified_cost_def}. Then, sequentially solve the LPs in \eqref{LP_start}-\eqref{LP_end} and \eqref{LP2_start}-\eqref{LP2_end}. Finally, synthesize the policy
$\pi^{\star}_{D,app}$$\in$$\Pi^{SD}(\mathcal{M})$ via the rule in \eqref{opt_policy_rule_approx}.

\subsubsection{Suboptimality analysis}
Let $c_{\min}$$:=$$\min\{c(s,a) : s\in S_r, a\in\mathcal{A}\}$ and $\widetilde{c}_{\max}$$:=$$\max\{\widetilde{c}(s,a) : s\in S_r, a\in \mathcal{A}\}$ be the minimum immediate cost and maximum modified immediate cost, respectively. Moreover, let $\underline{M}$ be a constant such that
\begin{align*}
    \underline{M}:= \min_{\pi\in \Xi^{SD}(\mathcal{M},B)} \mathbb{E}^{\pi}\Bigg[\sum_{t=1}^{\infty} f(s_t,a_t)\Bigg].
\end{align*}
$\underline{M}$ is the minimum number of steps (in expectation) taken under any deterministic policy before reaching the target set $B$ with maximum probability and can be computed efficiently via linear programming.

{\setlength{\parindent}{0cm}
\noindent \begin{prop}\label{approximation_thm} For any $\pi$$\in$$\Xi^{SD}(\mathcal{M},B)$, we have
\begin{align*}
  \underline{M}c_{\min} \overset{(1)}{\leq} J(\pi)  \overset{(2)}{\leq} \widetilde{J}(\pi) \overset{(3)}{\leq} M^{\star}\widetilde{c}_{\max}.
\end{align*}
\end{prop}}
\noindent \textbf{Proof:} The inequality (1) holds since, under the assumption that $c(s,a)$$=$$0$ for all $s$$\in$$B$$\cup$$S_0$, for any $\pi$$\in$$\Pi(\mathcal{M})$,
\begin{align*}
    J(\pi)\geq c_{\min} \mathbb{E}^{\pi}\Bigg[\sum_{t=1}^{\infty}\beta^{t-1}f(s_t,a_t)\Bigg] = \underline{M}c_{\min}.
\end{align*}
To show that the inequality (2) holds, note that
\begin{align*}
    J(\pi)=\sum_{t=1}^{\infty}\beta^{t-1}\text{Pr}^{\pi}\Big(S_t=s_t, A_t=a | S_1=s_1\Big) c(s_t,a_t).
\end{align*}
It follows from the definition of $T_{\min}(s)$ that we have $\text{Pr}^{\pi}(S_t=s_t, A_t=a | S_1=s_1)$$>$$0$ if and only if $t$$\geq$$T_{\min}(s)$. Then, for each nonzero term in the right hand side of the above equation, we have $\beta^{t-1}$$\leq$$\beta^{T_{\min}(s)-1}$ since $\beta$$<$$1$. Consequently, it follows from the definition of $\widetilde{c}$ that we have $J(\pi)$$\leq$$\widetilde{J}(\pi)$.

The inequality (3) holds for all $\pi$$\in$$\Xi^{SD}(\mathcal{M},B)$ since, under the assumption that $c(s,a)$$=$$0$ for all $s$$\in$$B\cup S_0$, we have
\begin{align*}
    \widetilde{J}(\pi)\leq \widetilde{c}_{\max} \mathbb{E}^{\pi}\Bigg[\sum_{t=1}^{\infty}f(s_t,a_t)\Bigg] \leq M^{\star}\widetilde{c}_{\max}. \qquad \Box
\end{align*}

{\setlength{\parindent}{0cm}
\noindent \begin{Corollary}
Let $\pi^{\star}_{D,app}$$\in$$\Pi^{SD}(\mathcal{M})$ be a policy satisfying the condition in \eqref{det_objective_approx}, and $\pi^{\star}_D$$\in$$\Pi^{SD}(\mathcal{M})$ be a policy satisfying the condition in \eqref{det_objective_2}. Then, we have
\begin{align}\label{corollary_ineq_1}
    J(\pi^{\star}_{D,app})-J(\pi^{\star}_D)\leq M^{\star} \widetilde{c}_{\max}- \underline{M} c_{\min}.
\end{align}
Furthermore, if the MDP $\mathcal{M}$ has only deterministic transitions, i.e., $\mathcal{P}_{s,a,s'}$$\in$$\{0,1\}$ for all $s$$\in$$S$ and for all $a$$\in$$\mathcal{A}$, then
\begin{align}\label{corollary_ineq_2}
    J(\pi^{\star}_{D,app})-J(\pi^{\star}_D)\leq \lvert S\rvert \widetilde{c}_{\max}.
\end{align}
\end{Corollary}}
\noindent \textbf{Proof:} It follows from Proposition \ref{approximation_thm} that $J(\pi)$$\leq$$\widetilde{J}(\pi)$$\leq$$M^{\star}\widetilde{c}_{\max}$ for all $\pi$$\in$$\Xi^{SD}(\mathcal{M},B)$. Then, we have $J(\pi^{\star}_{D,app})$$\leq$$M^{\star}\widetilde{c}_{\max}$. Similarly, since $\underline{M}c_{\min}$$\leq$$J(\pi)$ for all $\pi$$\in$$\Xi^{SD}(\mathcal{M},B)$, we have $\underline{M}c_{\min}$$\leq$$J(\pi^{\star}_D)$. Combining the above inequalities, we obtain the inequality in \eqref{corollary_ineq_1}. The result in \eqref{corollary_ineq_2} follows from the fact that $M^{\star}$$\leq$$\lvert S \rvert$ if the the MDP has only deterministic transitions. $\Box$

The bound in \eqref{corollary_ineq_2} shows that, for MDPs with deterministic transitions, the suboptimality of the proposed approximation algorithm grows \textit{at most} linearly in the size of the set $S$ of states.

\section{Simulation Results}

\newcommand{\StaticObstacle}[2]{ \fill[red] (#1+0.1,#2+0.1) rectangle (#1+0.9,#2+0.9);}
\newcommand{\initialstate}[2]{ \fill[black!30!brown] (#1+0.1,#2+0.1) rectangle (#1+0.9,#2+0.9);}
\newcommand{\goalstate}[2]{ \fill[black!50!green] (#1+0.1,#2+0.1) rectangle (#1+0.9,#2+0.9);}
\newcommand{\detectstatea}[2]{ \fill[blue!30] (#1+0.1,#2+0.1) rectangle (#1+0.9,#2+0.9);}
\newcommand{\detectstateb}[2]{ \fill[blue!20] (#1+0.1,#2+0.1) rectangle (#1+0.9,#2+0.9);}
\newcommand{\detectstatec}[2]{ \fill[blue!5] (#1+0.1,#2+0.1) rectangle (#1+0.9,#2+0.9);}
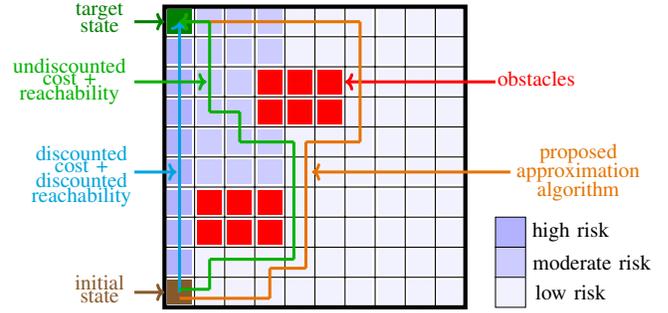
\begin{figure}[t!]
\centering
\scalebox{0.4}{
\begin{tikzpicture}
\draw[black,line width=1pt] (0,0) grid[step=1] (10,10);
\draw[black,line width=4pt] (0,0) rectangle (10,10);

\detectstatea{0}{1} \detectstatea{0}{2}
\detectstatea{0}{3} \detectstatea{0}{4}
\detectstatea{0}{5} \detectstatea{0}{6}
\detectstatea{0}{7} \detectstatea{0}{8}

\detectstatec{1}{0} \detectstatec{2}{0}
\detectstatec{3}{0} \detectstatec{4}{0}
\detectstatec{5}{0} \detectstatec{6}{0}
\detectstatec{7}{0} \detectstatec{8}{0}
\detectstatec{9}{0}

\detectstatec{1}{1} \detectstatec{2}{1}
\detectstatec{3}{1} \detectstatec{4}{1}
\detectstatec{5}{1} \detectstatec{6}{1}
\detectstatec{7}{1} \detectstatec{8}{1}
\detectstatec{9}{1}

\detectstatec{4}{2} \detectstatec{5}{2}
\detectstatec{6}{2} \detectstatec{7}{2}
\detectstatec{8}{2} \detectstatec{9}{2}
\detectstatec{4}{3} \detectstatec{5}{3}
\detectstatec{6}{3} \detectstatec{7}{3}
\detectstatec{8}{3} \detectstatec{9}{3}

\detectstatec{4}{4} \detectstatec{5}{4}
\detectstatec{6}{4} \detectstatec{7}{4}
\detectstatec{8}{4} \detectstatec{9}{4}
\detectstatec{4}{5} \detectstatec{5}{5}
\detectstatec{6}{5} \detectstatec{7}{5}
\detectstatec{8}{5} \detectstatec{9}{5}

\detectstatec{4}{8} \detectstatec{5}{8}
\detectstatec{6}{8} \detectstatec{7}{8}
\detectstatec{8}{8} \detectstatec{9}{8}
\detectstatec{4}{9} \detectstatec{5}{9}
\detectstatec{6}{9} \detectstatec{7}{9}
\detectstatec{8}{9} \detectstatec{9}{9}

\detectstatec{7}{6} \detectstatec{6}{6}
\detectstatec{9}{6} \detectstatec{8}{6}
\detectstatec{7}{7} \detectstatec{6}{7}
\detectstatec{9}{7} \detectstatec{8}{7}

\detectstateb{1}{9} \detectstateb{2}{9}
\detectstateb{3}{9} \detectstateb{1}{8}
\detectstateb{2}{8} \detectstateb{3}{8}
\detectstateb{1}{7} \detectstateb{2}{7}
\detectstateb{1}{6} \detectstateb{2}{6}
\detectstateb{1}{5} \detectstateb{2}{5}
\detectstateb{3}{5} \detectstateb{1}{4}
\detectstateb{2}{4} \detectstateb{3}{4}

\StaticObstacle{1}{2} \StaticObstacle{2}{2} \StaticObstacle{3}{2} \StaticObstacle{1}{3}
\StaticObstacle{2}{3} \StaticObstacle{3}{3}

\StaticObstacle{3}{6} \StaticObstacle{4}{6} \StaticObstacle{5}{6}
\StaticObstacle{3}{7} \StaticObstacle{4}{7} \StaticObstacle{5}{7}

\initialstate{0}{0}  \goalstate{0}{9}
			    
\draw[->, line width=1mm, color = black!30!brown] (-1,0.5) -- (0,0.5) node[below left = 0.5 cm and 2 cm, label={[align=center] \Huge initial\\ \Huge state}] {};

\draw[->, line width=1mm, color = black!50!green] (-1,9.5) -- (0,9.5) node[below left = 0.5 cm and 2 cm, label={[align=center] \Huge target\\ \Huge state}] {};

\draw[->, line width=1mm, color = red] (11,7.5) -- (6,7.5) node[below right = 0.3 cm and 6.2 cm, label={[align=center] \Huge obstacles}] {};

\draw[->, line width=1mm, color = black!10!cyan] (0.5,0.5) -- (0.5,9.5);

\draw[-, line width=1mm, color = black!10!orange] (0.5,0.3) -- (3.5,0.3);

\draw[-, line width=1mm, color = black!10!orange] (3.5,0.3) -- (3.5,1.3);

\draw[-, line width=1mm, color = black!10!orange] (3.5,1.3) -- (4.7,1.3);

\draw[-, line width=1mm, color = black!10!orange] (4.7,1.3) -- (4.7,5.5);

\draw[-, line width=1mm, color = black!10!orange] (4.7,5.5) -- (6.5,5.5);

\draw[-, line width=1mm, color = black!10!orange] (6.5,5.5) -- (6.5,9.5);

\draw[->, line width=1mm, color = black!10!orange] (6.5,9.5) -- (0.5,9.5);

\draw[-, line width=1mm, color = black!30!green] (0.5,0.6) -- (1.5,0.6);

\draw[-, line width=1mm, color = black!30!green] (1.5,0.6) -- (1.5,1.6);

\draw[-, line width=1mm, color = black!30!green] (1.5,1.6) -- (4.3,1.6);

\draw[-, line width=1mm, color = black!30!green] (4.3,1.6) -- (4.3,5.5);

\draw[-, line width=1mm, color = black!30!green] (4.3,5.5) -- (2.5,5.5);

\draw[-, line width=1mm, color = black!30!green] (2.5,5.5) -- (2.5,6.5);

\draw[-, line width=1mm, color = black!30!green] (2.5,6.5) -- (1.5,6.5);

\draw[-, line width=1mm, color = black!30!green] (1.5,6.5) -- (1.5,9.5);

\draw[->, line width=1mm, color = black!30!green] (1.5,9.5) -- (0.5,9.5);

\draw[fill = blue!30] (11,2) rectangle (12,3);
\draw[fill = blue!20] (11,1) rectangle (12,2);
\draw[fill = blue!5] (11,0) rectangle (12,1);

\node at (13.5,2.5) {\Huge high risk};
\node at (14.2,1.5) {\Huge moderate risk};
\node at (13.5,0.5) {\Huge low risk};

\draw[->, line width=1mm, color = black!10!orange] (11.5,4.5) -- (4.9,4.5) node[below right = 1.1 cm and 8.7 cm, label={[align=center] \Huge proposed\\ \Huge approximation \\ \Huge algorithm}] {};

\draw[->, line width=1mm, color = black!10!cyan] (-1,4.5) -- (0.4,4.5) node[below left = 1.3 cm and 3 cm, label={[align=center] \Huge discounted \\ \Huge cost +\\ \Huge discounted \\ \Huge reachability}] {};

\draw[->, line width=1mm, color = black!30!green] (-1,7.5) -- (1.4,7.5) node[below left = 0.9 cm and 4.4 cm, label={[align=center] \Huge undiscounted \\ \Huge cost +\\ \Huge reachability}] {};

\end{tikzpicture}

}
\caption{An illustration of the environment in the motion planning example. The agent aims to reach the target state while avoiding risky regions for as long as possible. By combining a discounted cost criterion with a reachability criterion, the proposed approximation algorithm generates a trajectory that postpones the entrance to risky regions while ensuring the satisfaction of the task. }
\label{simulation_tikz}
\end{figure}
In this section, we showcase the performance of the proposed algorithms on a motion planning example. We run the computations on a 3.2 GHz desktop with 8 GB RAM and utilize the GUROBI solver \cite{gurobi} for optimization.

Consider an autonomous vehicle that aims to deliver a package to a target region in an adversarial environment. Specifically, in each region in the environment, the vehicle faces the risk of being attacked by an adversary whose objective is to prevent the vehicle from reaching its target. Being aware of the threat, the vehicle's objective is to reach its target by following minimum-risk trajectories.

We model the environment as a grid world illustrated in Figure \ref{simulation_tikz}. The agent starts from the initial state (brown) and aims to reach the target state (green) while avoiding the states that are occupied by obstacles (red). The set of actions available to the agent are $\{up,down, left,right, stay\}$. Under an action $a$$\in$$\{up,down, left,right\}$, the agent moves to the successor state in the desired direction with probability 0.9 and stays in its current state with probability 0.1. Under the $stay$ action, the agent stays in its current state with probability 1. Note that, due to the stochasticity in action outcomes, the agent
has an infinite decision horizon for reaching the target state with probability 1.

Each state in the environment has an associated cost that represents the risk of being attacked by an adversary. The states along the minimum length trajectory have a high risk, i.e., $c(s,a)$$=$$4$. The states that are between the target and the obstacles have a moderate risk, i.e., $c(s,a)$$=$$2$. Finally, all other states have a low risk, i.e., $c(s,a)$$=$$1$.

We synthesize three stationary deterministic policies for the agent using two existing methods in the literature and the proposed approximation algorithm. With 309 binary variables in the MILP formulation, the computation exceeds the memory limit after 1325 seconds.

We first synthesize a policy based on the classical constrained MDP formulation in which both the incurred costs and the collected rewards are discounted \cite{altman1999constrained}. In particular, we express the reachability criterion as a total reward criterion as shown in \eqref{equivalence_reach}. We then discount both the costs and the rewards with $\beta$$<$$1$ and synthesize a policy that minimizes the total discounted cost among the ones that maximizes the total discounted reward. The agent's trajectory under the synthesized policy is shown in Figure \ref{simulation_tikz} (as ``discounted cost + discounted reachability"). Due to the discounting in the rewards, the agent aims to reach the target as quickly as possible. In other words, the discounting results in a behavior that assigns an artificial importance to the reachability objective, which generates a trajectory that visits the high-risk regions in the environment.

The second policy is based on a multi-objective MDP formulation \cite{chatterjee2006markov} in which neither the costs nor the rewards are discounted. Specifically, we express the reachability criterion as a total reward criterion and consider undiscounted costs by setting $\beta$$=$$1$. The agent's trajectory under the synthesized policy is shown in Figure \ref{simulation_tikz} (as ``undiscounted cost + reachability"). The agent minimizes the total risk along its trajectory. However, due to lack of discounting in the costs, it cannot adjust the time it enters the states with moderate risk. As a result, the synthesized policy generates a trajectory that visits more states with moderate risk than the minimum achievable one.

We synthesize the third policy through the approximation algorithm presented in Section \ref{app_algo_sec} by choosing $\beta$$=$$0.9$. The agent's trajectory under the synthesized policy is shown in Figure \ref{simulation_tikz} (as ``proposed approximation algorithm"). The synthesized policy attains the same total risk with the second policy explained above. However, thanks to the discounting in the costs, the agent establishes a hierarchy between the policies with the same total risk and follows a trajectory that visits the minimum number of states with moderate risk.

The above example illustrates the two main benefits of the formulation that combines a discounted cost criterion with a probabilistic reachability criterion. First, by not discounting the reachability constraint, the agent is able to follow long trajectories if it is optimal to do so. Second, by discounting the costs, the agent is able to create an ordering between the policies that incur the same total cost.

\section{Conclusions and Discussions}
We studied the problem of synthesizing a policy in a Markov decision process (MDP) 
which, when followed, causes an agent to
reach a target state in the MDP while minimizing its total discounted cost. We showed that, in general, an optimal policy for this problem might not exist, but there always exists a near-optimal policy which can be synthesized efficiently. We also considered the synthesis of an optimal stationary deterministic policy and established that the synthesis of such a policy is NP-hard. Finally, we proposed a linear programming-based algorithm to synthesize stationary deterministic policies with a theoretical suboptimality guarantee.

In the second part of the paper, we presented an analysis over \textit{stationary} deterministic policies. In general, \textit{Markovian} deterministic policies achieve lower total discounted costs than stationary deterministic policies when subject to a reachability constraint. For such instances, one can synthesize Markovian policies by taking a Cartesian product of the MDP with a deterministic finite automaton representing the ``time-dependency'' of the policy and applying the algorithms developed in this paper.

We have no proof showing that the suboptimality bound for the proposed approximation algorithm is tight. Considering the generality of the formulation studied, it may be of interest to construct an MDP instance to prove the tightness of the derived bound or to develop approximation algorithms with better performance guarantees.

  \bibliographystyle{IEEEtran}
\bibliography{main.bib}
\end{document}